\documentclass[10pt]{amsart}

\usepackage{latexsym}
\usepackage{amssymb}
\usepackage{graphicx}
\usepackage{latexsym}
\usepackage{amssymb}

\newcommand{\excise}[1]{}

\numberwithin{equation}{section}

\newtheorem{thm}{Theorem}[section]
\newtheorem{lemma}[thm]{Lemma}

\newtheorem{prop}[thm]{Proposition}

\newtheorem{theorem}{Theorem}
\newtheorem{Example}[thm]{Example}
\newtheorem{Remark}[thm]{Remark}
\newtheorem{Alg}[thm]{Algorithm}
\newtheorem{Defn}[thm]{Definition}

\newenvironment{remark}{\begin{Remark}\rm}
                {\mbox{}~\hfill$\square$\end{Remark}}

\newenvironment{defn}{\begin{Defn}\rm}
        {%\mbox{}~\hfill$\square$
         \end{Defn}}

    {\begin{eqnarray}}
    {\end{eqnarray}$\!\!$}

\newenvironment{eq*}%
    {\begin{eqnarray*}}
    {\end{eqnarray*}$\!\!$}

    {\begin{equation}}
    {\end{equation}$\!\!$}

\newenvironment{eqn*}%
    {\begin{equation*}}
    {\end{equation*}$\!\!$}

        {\begin{list}
                {\noindent\makebox[0mm][r]{\arabic{enumi}.}}
                {\leftmargin=5.5ex \usecounter{enumi}}
        }
        {\end{list}}

        {\begin{list}
                {\noindent\makebox[0mm][r]{(\roman{enumi})}}
                {\leftmargin=5.5ex \usecounter{enumi}}
        }
        {\end{list}}

\def\<{\langle}
\def\>{\rangle}
\def\0{{\mathbf 0}}
\def\1{{\mathbf 1}}

\def\CC{{\mathbb C}}

\def\RR{{\mathbb R}}

\def\into{\hookrightarrow}

\def\lieT{{\mathfrak t^*}}
\def\lieS{{\mathfrak s^*}}
\def\red{{/\!/\!}}
\def\cone{{\Sigma}}
\def\dualcone{{\Lambda}}
\def\pr{{\rm \rho}}

\begin{document}%%%%%%%%%%%%%%%%%%%%%%%%%%%%%%%%%%%%%%%%%%%%%%%%%%%%%%%%%
%%%%%%%%%%%%%%%%%%%%%%%%%%%%%%%%%%%%%%%%%%%%%%%%%%%%%%%%%%%%%%%%%%%%%%%%%

\title[Localization theorems by symplectic cuts]%
    {Localization theorems by symplectic cuts}
\author{Lisa Jeffrey}
\thanks{LJ was supported by a grant from NSERC.
This work was begun during a  visit of LJ to Harvard
University during the
spring of 2003, whose support during this period is
acknowledged.}
\address{University of Toronto\\Toronto, ON\\Canada}
\email{jeffrey@math.toronto.edu}

\author{Mikhail Kogan}
\thanks{MK was supported by the National Science Foundation under
agreement No. DMS-0111298}
\address{Institute for Advanced Study\\Princeton, NJ\\USA}
\email{mish@ias.edu}

\begin{abstract}
 Given a compact symplectic manifold $M$ with the
Hamiltonian action of a  torus $T$,  let zero be a regular value
of the moment map, and $M_0$ the symplectic reduction at zero.
Denote by~$\kappa_0$ the Kirwan map $H^*_T(M)\to H^*(M_0)$.  For an
equivariant cohomology  class $\eta\in H^*_T(M)$ we present new
localization formulas which express $\int_{M_0} \kappa_0(\eta)$ as
sums of certain integrals over the connected components of the
fixed point set $M^T$. To produce such a formula we apply a
residue operation to the Atiyah-Bott-Berline-Vergne
 localization formula for an equivariant
form on the symplectic cut of $M$ with respect to a certain cone, and
then, if necessary, iterate this process using other cones. When
all cones used to produce the formula are one-dimensional we
recover, as a special case, the localization formula of Guillemin
and Kalkman~\cite{GK:residue}. Using similar ideas, for a special
choice of the cone (whose dimension is equal to that of $T$) we
give a new proof of the Jeffrey-Kirwan localization
formula~\cite{JK:LocalizationNonabelian}.
\end{abstract}

\maketitle

\renewcommand{\t}{{\mathfrak t}}
\newcommand{\tc}{{\mathfrak t_\CC}}
\newcommand{\res}{{\rm Res}}
\newcommand{\liet}{{\bf t}}
\newcommand{\liets}{{\bf t^*}}
\newcommand{\C}{ {\bf C}}
\newcommand{\R}{ { \bf R}}
\newcommand{\Jac}{{\bigtriangleup}}

{\center{ \em  This paper is dedicated to Alan Weinstein
on the occasion of his 60th birthday.}}

\section{Introduction}
\label{sec:intro}

Assume we are given a compact symplectic manifold $M$ with the
Hamiltonian action of a torus~$T$. There are two kinds of
localization theorems which express the integral over $M$ of an
equivariant cohomology class $\eta \in H^*_T(M)$ and the integral
over the reduced space $M_0$ of the Kirwan map $\kappa_0(\eta)$ as
sums of certain terms which involve integration over the connected
components of the fixed point set $M^T$. In particular, these
localization theorems say that both $\int_M\eta$ and
$\int_{M_0}\kappa_0(\eta)$ depend only on the restriction of
$\eta$ to $M^T$.

More specifically, the localization theorem of
Atiyah-Bott~\cite{AB} and  Berline-Vergne~\cite{BV}
(or just  the ABBV localization theorem) expresses the
integral over $M$ (that is, the pushforward in equivariant cohomology with
respect to the map $M\to {\rm pt}$) of a class
$\eta\in H^*_T(M,\CC)$ as a sum of integrals over the connected components
$F$ of the fixed point set $M^T$ of the restriction of $\eta$ to $F$
divided by the
equivariant Euler class of the normal bundle of $F$:
$$
\int_M\eta= \sum_{F} \int_{F}
\frac{\iota_{F}^*(\eta)}{e(\nu(F))},
$$
where $\iota_F$ is the natural inclusion $F\to M$ and $e(\nu(F))$ is the
equivariant Euler class of the normal bundle $\nu(F)$.

To treat the other kind of localization theorems, let zero be a
regular value of the moment map $\mu:M\to \mathfrak t$,
$(M_0,\omega_0)$ the symplectic reduction at zero and
$\kappa_0:H^*_T(M)\to H^*(M_0)$ the Kirwan map. Then the
Jeffrey-Kirwan~\cite{JK:LocalizationNonabelian} and
Guillemin-Kalkman~\cite{GK:residue} localization theorems express
the integral over $M_0$ of classes $\kappa_0(\eta)e^{\omega_0}$
and $\kappa_0(\eta)$ respectively as sums over some connected
components of $M^T$ of certain terms similar to those appearing in
the ABBV localization theorem. We postpone the precise statement
of the Guillemin-Kalkman theorem until Section~\ref{sec:GK}.
However the Jeffrey-Kirwan localization theorem applied to the
case of abelian group actions states the following under the above
assumptions.

\begin{theorem}\label{thm:main}{\rm \cite{JK:LocalizationNonabelian}}
For $\eta\in H^*_T(M)$ we have
\begin{equation}
\label{eq:main} \int_{M_0}\kappa_0(\eta)e^{\omega_0} = c\cdot
\res^\Lambda \left(\sum_{F} e^{i(\mu(F))(X)}\int_F\frac{\iota^*_F
(\eta(X) e^{\omega})}{e(\nu(F))(X)}[dX]\right)
\end{equation}
where $c$ is a non-zero constant, $X$ is a variable in
$\mathfrak{t}\otimes \mathbb C$ so that a $T$-equivariant
cohomology class can be evaluated at $X$, and $\res^\Lambda$
is a multi-dimensional
residue with respect to the cone $\Lambda\subset\t$ defined in
Section~\ref{sec:residues}.
\end{theorem}

The similarity of the ABBV localization theorem to the
Jeffrey-Kirwan localization theorems is transparent, and is  of
course not a coincidence. In the case of Hamiltonian circle
actions, Lerman~\cite{Lerman} showed that it is possible to deduce
the Jeffrey-Kirwan and Guillemin-Kalkman theorems from ABBV
localization using the techniques of symplectic cutting. The idea,
which was also present, though not explicitly stated,
in~\cite{GK:residue}, is to use symplectic cutting to produce a
Hamiltonian space whose connected components of the $T$-fixed set
are either the reduced space $M_0$ or some connected components of the
set $M^T$. Then the ABBV localization theorem on this symplectic
cut yields a formula which relates integration over some of the
connected components of $M^T$ and over the reduced space $M_0$. To
arrive at the Jeffrey-Kirwan and Guillemin-Kalkman theorems for
the case of a circle action, it remains to apply residues to both
sides of the ABBV localization formula for the symplectic cut.

The goal of this paper is to illustrate that an analogous approach
works for higher dimensional torus actions as well. Let us outline
the main idea. Lerman's original definition of
%\comment{Misha --
%do you mean Lerman's definition of symplectic cutting or his proof
%of GK using symplectic cutting? Lisa, In this sentence I only mean
%Lerman's definition of symplectic cutting.}
symplectic cutting was
given for the case of circle actions. However it can be
generalized to multidimensional torus actions, when the symplectic
cut is defined using any rational convex polytope. We will only
consider the special case of symplectic cutting with respect to a
cone. Let $\cone$ be a convex rational polyhedral cone (centered
at the origin) in the dual $\mathfrak t^*$ of the Lie algebra
$\mathfrak t$. If $\sigma$ is an open face of $\cone$ (that is,
the interior of a face of $\cone$), let $T^\sigma$ be the subtorus
of $T$ whose Lie algebra is annihilated by $\sigma$. Then, as a
topological space, the symplectic cut $M_\cone$ is
$\mu^{-1}(\cone)/\sim$, where $p\sim q$ if $\mu(p)=\mu(q)\in
\sigma$, for some open face $\sigma$ of $\cone$ and $p,q$ lying in
the same $T^\sigma$ orbit. As shown in~\cite{LMTW}, for a generic
choice of~$\cone$, the cut space $M_\cone$ is a symplectic
orbifold with a Hamiltonian $T$ action. The moment map image of
$M_\cone$ is just the intersection $\mu(M)\cap \cone$. Moreover,
any equivariant cohomology class $\eta\in H^*_T(M)$ naturally
descends to an equivariant class $\eta_\cone$ on $M_\cone$.

Some connected components of the fixed point set $M_\cone^T$ may
be identical to those of $M^T$; we call them the {\it old}
connected components (see Definition~\ref{def:points}). One of
the connected components of~$M_\cone^T$ is always the reduced
space $M_0$. Hence if we apply ABBV localization to the class
$\eta_\cone$ on the symplectic cut $M_\cone$, we will get a
formula which relates the integration over $M_0$ to the integration
over the connected components of $M_\cone^T$, some of which are
the connected components of $M^T$. We will show that we can apply
the iterated residue operation to both sides of this ABBV
localization formula, so that
the term corresponding to $M_0$ simplifies. More specifically,
this term  becomes a constant times the
integral of $\kappa_0(\eta)$ over $M_0$.

The major difference with the circle case is that in the formula
just described, besides a contribution from the term corresponding
to $M_0$ and the terms corresponding to the old connected
components, there will be contributions coming from the {\it new}
connected components of $M^T_\cone$, which are neither part of
$M^T$ nor part of $M_0$ (see Definition~\ref{def:points}). However, using
the ideas of~\cite{GK:residue} we can iterate this process to get
rid of these terms. Namely for each new connected component $F'$,
we can symplectically cut a certain submanifold $M'$ of $M$ with
respect to some cone $\cone'$, so that~$F'$ is a connected
component of the $T$ action on $M'_{\cone'}$. Then we can apply
the ABBV theorem and the residue operation again to $M'_{\cone'}$
to express the integral over~$F'$ in terms of integrals over the
connected components of the fixed point set $M'^T_{\cone'}$. As we
will show this process can be iterated until all the terms coming
from the new fixed points disappear. So the integration of the
Kirwan map over~$M_0$ can be expressed as a sum of terms which involve
integration only over the connected components of~$M^T$.

The objects which carry the information about which cones are
chosen in this process are called {\em dendrites}, and are defined in
Section~\ref{sec:dendrites}. Every dendrite gives a
localization formula. If all the cones used in a dendrite are
one-dimensional, then we recover the
Guillemin-Kalkman~\cite{GK:residue} localization theorem. However,
if we choose higher dimensional cones, in other words multi-dimensional
dendrites, we get new localization formulas.

Notice  that the Jeffrey-Kirwan localization theorem  does not involve any
iteration. Nevertheless, it fits into the  framework just
described. While for actions of tori of dimension greater than one
any symplectic cut with respect to a cone always has new connected
components of $M_\cone^T$  (the connected components which are
neither in $M^T$ nor in $M_0$), it is plausible that their
contribution to the ABBV localization theorem becomes zero after
taking residues. In this case the iterative process described
above would stop at step one. We were not able to use precisely this
argument to give a proof of the Jeffrey-Kirwan localization theorem.
However, we
will show that for a good choice of the cone~$\cone$, a very
similar argument which involves taking residues of the ABBV formula
for symplectic cuts yields the Jeffrey-Kirwan formula given in
Theorem~\ref{thm:main}.

The paper is organized as follows. In Section~\ref{sec:cuts} we
carefully review  well-known objects from the theory of
Hamiltonian group actions, such as symplectic reductions,
symplectic cuts, equivariant cohomology and the  Kirwan map. We
also  recall the orbifold version of the ABBV localization
theorem. Section~\ref{sec:residues} is devoted to the residue
operation. The definitions of  residues are based on the theory
of complex variables, do not involve any symplectic geometry and
are independent of the material summarized
in  Section~\ref{sec:cuts}. In Section~\ref{sec:GK}
we present the generalization of the Guillemin-Kalkman
localization theorem to the case of multi-dimensional dendrites.
Finally, using an analogous approach, in Section~\ref{sec:JK} we give
a new proof of the Jeffrey-Kirwan localization theorem.

\section{Symplectic cuts and other preliminaries}
\label{sec:cuts} In this section we recall the construction of
symplectic cuts with respect to cones and other results in
symplectic geometry. All of the results in this section (with only
one exception: Proposition~\ref{prop:euler}) have appeared in the
literature, so we state them without proof.

\subsection{Symplectic reduction.}

Let $(M,\omega)$ be a symplectic manifold with  the
action of a Hamiltonian torus
$T$ and  a moment map $\mu: M\to \lieT$. For $p\in \lieT$
\emph{the symplectic reduction $M_p=M\red_p T$ of $M$ at $p$} is
defined to be $\mu^{-1}(p)/T$. Whenever $p$ is a regular value of
the moment map, the symplectic reduction $M_p$ is an orbifold.

For a subtorus $H\subseteq T$, denote by $M^H$ the fixed point set
of the $H$ action on $M$ and by $M^H_i$ the connected components
of $M^H$. It is well known that every $M^H_i$ is a symplectic
manifold with a Hamiltonian action of the torus $T/H$. The convexity
theorem of
Atiyah~\cite{A} and  Guillemin-Sternberg~\cite{GS}
states that if $M$ is compact and connected then $\mu(M)$ is a
convex polytope. In particular, every~$\mu(M^H_i)$ is a convex
polytope inside $\mu(M)$;  we call it \emph{a wall of $\mu$}.

Let $T$ be a product of two subtori $H\times S$. Denote by $\pi$
the natural projection $\lieT\to \lieS$. Then the composition
$\pi\circ \mu$ is a moment map for the $S$ action on $M$.
Moreover, the action of~$S$ on $M$ restricts to a Hamiltonian
action on  $M^H$ whose moment map is again given by
$\pi\circ\mu$. Because of this, for $p\in \mu(M^H)$ we call the
space
$$
M^H\red_p S= (\mu^{-1}(p)\cap M^H)\big/ T=(\mu^{-1}(p)\cap
M^H)\big/ S
$$
the symplectic reduction of $M^H$ at $p$.

If $M$ is compact and $q\in \lieS$ then the set $\pi^{-1}(q) \cap
\mu(M^H)$ contains finitely many points $p_i$. It is easy to see
that the fixed point set of the $H$ action on $M\red_q S$ is the
union of spaces $M^H\red_{p_i}S$.

\subsection{Cutting with respect to a cone}
\label{subsec:cut} Given a linearly independent set $\beta
=\{\beta_1,\dots,\beta_k\}$ of weights of~$T$, define the cone
$\cone=\cone_\beta\subset \lieT$ to be given by all nonnegative
linear combinations of the weights of $\beta$, namely
$\cone=\{\sum s_i\beta_i| s_i\geq0\}$. (Note that we allow $k$ to be less
than $\dim T$.) For a set~$I$ of indices
between $1$ and $k$, denote by $\cone^I$ the open face of $\cone$
given by positive linear combinations of weights indexed by the
elements of $I$, that is, $\cone^I=\{\sum_{i\in I} s_i\beta_i| s_i>
0\}$. It will be convenient to denote the subtorus of $T$
perpendicular to $\cone^I$ by~$T^I$.

Assume that $(M,\omega)$ carries a Hamiltonian $T$ action with
moment map~$\mu$. For simplicity assume that the $T$ action is effective
so that $\mu(M)$ is a polytope of dimension equal to $\dim T$. As a
topological space the symplectic cut with respect to a cone
$\cone=\cone_\beta$ is the space
$$
M_\cone =\mu^{-1}(\Sigma)/\sim,
$$
where $x\sim x'$ if $\mu(x)=\mu(x')\in \cone_\beta^I$ and $x$ and
$x'$ lie in the same $T^I$ orbit. Clearly, the torus action on~$M$
descends to an action on~$M_\cone$. Notice that the action on
$M_\cone$ is not  effective  unless $\dim \cone=\dim T$,
since the subtorus $T^\cone$ of~$T$ perpendicular to $\cone$ acts
trivially on $M_\cone$.

In the case of a circle action, when every cone is just a ray,
Lerman \cite{Lerman} realized that under certain mild conditions
this space is an orbifold, it carries a natural symplectic form,
and the residual torus action is Hamiltonian. As shown in
\cite{LMTW} similar results hold for  symplectic cuts with
respect to cones. Even more generally, symplectic cutting has been
extended to cutting with respect to polytopes \cite{LMTW}, cutting
for nonabelian groups~\cite{Meinrenken} and K\"{a}hler
cutting~\cite{BGL}. For the purposes of this paper we only need to
consider symplectic cutting with respect to cones.

To recall the results of~\cite{LMTW} let us give another
construction of $M_\cone$. Let $\CC_{\beta_i}$ be a complex line
on which $T$ acts with weight $\beta_i$. Let
$\CC_\beta=\CC_{\beta_1}\times\dots\times\CC_{\beta_k}$. Then the
$T$ action on $\CC_\beta$ is Hamiltonian with respect to the
symplectic form $\omega_\CC=\sqrt{-1}\sum dz_i\wedge d\bar z_i$,
where $z_i$ is the standard complex coordinate on $\CC_{\beta_i}$.
Its moment map is $\psi(z)=\sum \beta_i|z_i|^2$, so that the image
of $\CC_\beta$ under the moment map is the cone $\cone$. Consider
the symplectic form $(\omega,-\omega_\CC)$ on $M\times \CC_\beta$.
Then  the diagonal torus $\Delta T\subset T\times T$ acts on this
space in a Hamiltonian fashion and its moment map $\phi$ is given
by $\phi(x,z) =\mu(x)-\psi(z)$. As shown in~\cite{LMTW} the
symplectic cut $M_\cone$ is homeomorphic to the symplectic
reduction $(M\times \CC_\beta)\red_0 \Delta T$.

Hence, to guarantee that $M_\cone$ is a symplectic orbifold it is
enough to assume that zero is a regular value of $\phi$. It is
easy to see that this is equivalent to the following
\renewcommand{\labelitemi}{($\star$)}
\begin{itemize} \item
Every $\cone^I$ is transverse to every wall $W$ of $\mu(M)$, that is,
 $\dim(\cone^I\cap W)\leq |I| +\dim W -\dim T$.
\end{itemize}
If ($\star$) holds we say that $\cone$ is {\emph{transverse to
$\mu$} and from now on we consider only the cones satisfying~($\star$).

We can write $T\times T$ as the product of the first copy of $T$
(that is, $T\times e$) and $\Delta T$. Hence after reducing
$M\times \CC_\beta$ with respect to~$\Delta T$, the action of the
first copy of $T$ descends to an action on the
reduction~$M_\cone$. As shown in~\cite{LMTW}, this action is
Hamiltonian and there exists a moment map~$\mu_{\cone}$ on
$M_\cone$ whose moment map image is the intersection $\mu(M)\cap
\cone$.

Let us describe the connected components of the fixed point set
$M^T_\cone$. We will separate them into three sets: the old fixed
points, the new fixed points, and the fixed points at zero.
\begin{defn}
\label{def:points} The old fixed points exist only when $\dim
\cone=\dim T$; they are all the connected components~$F_i$
of~$M^T_\cone$ for which $\mu_{\cone}(F_i)$ is in the interior of
$\cone$. The set of fixed points at zero is defined to be the
connected component~$\mu_{\cone}^{-1}(0)$ of~$M^T_\cone$, which is
just the symplectic reduction $M_0$.  The new fixed points are all
the other connected components $F_i'$ of $M^T_\cone$.
\end{defn}
It is not difficult to see that every $F_i$ is also a connected
component of the fixed point set $M^T$,
while $F_i'$ do not
correspond to any fixed points on $M$.

\subsection{Kirwan map}
We recall the definition of  equivariant cohomology $H^*_T(M)$
with complex coefficients using the Cartan model. Denote
$$
\Omega^*_T(M)=S(\lieT)\otimes\Omega^*(M)^T,
$$
where $S(\lieT)$ denotes the algebra of polynomials on $\mathfrak
t$ and $\Omega^*(M)^T$ denotes all $T$ invariant differential forms.
So, if $f\in S(\lieT)$, $\alpha\in\Omega^*(M)^T$  and $X\in\mathfrak
g$ we set $(f\otimes \alpha)(X) =f(X)\alpha$. The $T$-equivariant
differential on $\Omega^*_T(M)$ which defines the equivariant
cohomology is given by
$$
d_T(f\otimes \alpha)  = f\otimes d\alpha - \sum x_i
f\otimes\imath_{\xi^i_M} \alpha
$$
where $\xi_1,\dots,\xi_d$ is a basis of $\mathfrak g$; $x_1,\dots,
x_d$ is the dual basis of $\lieT$;  and $\xi^i_M$ is the vector
field generated by the action of $\xi^i$. If the equivariant form
$\beta\in\Omega^*_T(M)$ is closed, that is $d_T\beta=0$, we denote
by $[\beta]$ the cohomology class it represents.

If $T$ acts locally freely on $M$ then it is well known that
$$
H^*_T(M)=H^*(M/T),
$$
and if the action is trivial then
$$
H^*_T(M)=H^*(M)\otimes S(\lieT) =H^*(M)\otimes H^*_T(pt).
$$

For $p\in \lieT$ let $i_p$ be the inclusion $\mu^{-1}(p)\into M$.
If $p$ is a regular value of the moment map then $T$ acts locally
freely on $\mu^{-1}(p)$ so that
\begin{equation}
\label{eq:isom} H^*_T(\mu^{-1}(p))\cong H^*(\mu^{-1}(p)/T)=
H^*(M\red_pT)
\end{equation}
Composition of the pullback $i^*_p$ with the isomorphism
(\ref{eq:isom}) defines the Kirwan map
$$
\kappa_p: H^*_T(M)\to H^*(M\red_pT).
$$
Kirwan~\cite{Kirwan} showed that if $M$ is compact this map is
surjective. In the presence of another Hamiltonian action on $M$
by a torus $T'$ which commutes with the action of $T$, the Kirwan
map generalizes to its equivariant version:
$$
\kappa_p: H^*_{T\times T'}(M)\to H^*_{T'}(M\red_pT).
$$
Kirwan surjectivity was generalized to this case in~\cite{Goldin}.

Analogously, in the case when $T=H\times S$ and $p\in \mu(M^H)$  we can
define the Kirwan map
$$
\kappa_p^H: H^*_S(M^H)\to H^*(M^H\red_p S).
$$
We will mostly be interested in the equivariant version of this
map. Let us take  account of the trivial action of
$H$ on both $M^H$ and $M^H\red_p S$. So in the rest of the paper
$\kappa_p^H$ will be the map
\begin{equation}
\label{eq:kappa_H} \kappa_p^H: H^*_T(M^H)\to H^*_H(M^H\red_p S).
\end{equation}

Let us now apply the equivariant version of the Kirwan map to
symplectic cuts with respect to cones. Given a cone
$\cone=\cone_\beta$, the product $T\times T$ acts on $M\times
\CC_\beta$ and the symplectic cut $M_\cone$ is produced by
reducing $M\times \CC_\beta$ at~$0$. So if we think of $T\times
T$ as the product of $T\times e$ and $\Delta T$, then the
equivariant version of the Kirwan map produces the map
$$
\kappa_\cone: H^*_{T\times T}(M\times \CC_\beta)\to H^*_{T\times
e} (M_\cone)=H^*_{T} (M_\cone).
$$
As mentioned before the action of $T$ on $M_\cone$ might not be
effective, since the torus $T^\cone$ orthogonal to the cone
$\cone$ acts trivially on $M_\cone$. In particular,
$$
H^*_T(M_\cone)=H^*_{T/T^\cone}(M_\cone)\otimes H^*_{T^\cone}(pt).
$$
For $\eta\in H^*_T(M)$ denote by $\eta_\cone$ the class
$\kappa_\cone (\eta\otimes 1)\in H^*_T(M_\cone)$. It is an easy
exercise to see that $\eta_\cone\in
H^*_{T/T^\cone}(M_\cone)\otimes 1$.

Let us also notice that if $F_i$ is a connected component for the
fixed point set of both $M$ and $M_\cone$ with
$\mu(F_i)=\mu_{\cone}(F_i)\in {\rm Int}(\cone)$, then
$$
\iota^*_{F_i}\eta = \iota^*_{F_i} \eta_\cone
$$
where by abuse of notation we denote by $\iota_{F_i}$ both the
inclusion of $F_i$ into $M$ and that into $M_\cone$. (Because $\cone$
is transverse to $\mu(M)$, the existence of such $F_i$ implies
that the dimension of the cone is at least the dimension of
$\mu(M)$.) It is also easy to see that
$$
\iota_{M_0}\eta_\cone=\kappa_0(\eta)\otimes 1 \in H^*(M_0)\otimes
1\subset H^*_T(M_0),
$$
where, since $T$ acts trivially on $M_0$, we know that
$H^*_T(M_0)=H^*(M_0)\otimes S(\lieT)$.

\subsection{ABBV localization theorem}
Suppose that $M$ is compact, oriented and carries the action of a
torus $T$. The pushforward map from $H^*_T(M)$ to $H^*_T(pt)$ is
the integration over~$M$ and is denoted by~$\int$. The
theorem of Atiyah-Bott
and Berline-Vergne (or just the ABBV theorem)~\cite{AB,BV}
states that for~$\eta\in H^*_T(M)$
\begin{equation}
\label{eq:ABBV} \int_M\eta= \sum_{F_i} \int_{F_i}
\frac{\iota_{F_i}^*(\eta)}{e(\nu(F_i))}.
\end{equation}
Here $F_i$ are the connected components of the fixed point set
$M^T$, $\iota_{F_i}:F_i\into M$ are their inclusions, $\nu(F_i)$
are the normal bundles of $F_i$, and $e(\nu(F_i))$ are their
$T$-equivariant Euler classes.

We need to know more about these $T$-equivariant Euler
classes. Because of the splitting principle~\cite{BT} we can
assume without loss of generality that $e(\nu(F_i))$ splits as a
sum of line bundles $L_1\otimes\dots\otimes L_k$. Assume $T$ acts
on the fibers of $L_i$ with weight $\lambda_i$. Then the
$T$-equivariant Euler class~is
\begin{equation}
\label{eq:Euler-class} e(\nu(F_i))=\prod(\lambda_i + c_1(L_i)),
\end{equation}
where $c_1(L_i)$ is the first Chern class of $L_i$.

The ABBV localization formula was generalized to orbifolds by
Meinrenken~\cite{Meinrenken}. We refer to~\cite{Meinrenken} for
details. Let us just mention that the only difference
with~(\ref{eq:ABBV}) is the appearance of constants before each
term of the formula:
\begin{equation}
\label{eq:ABBV-orbifold} \frac{1}{d_M}\int_M\eta= \sum_{F_i}
\frac{1}{d_{F_i}} \int_{F_i}
\frac{\iota_{F_i}^*(\eta)}{e(\nu(F_i))},
\end{equation}
where for a connected orbifold $X$, the size of the finite
stabilizer at a generic point of $X$ is equal to~$d_X$.
Moreover~(\ref{eq:Euler-class}) is still a valid formula for
orbifolds, where  the $\lambda_i$ (when properly interpreted) are
rational weights of the action on the normal bundle.

Another important property of the equivariant Euler classes is the
following generalization of \cite[Proposition~3.1]{GK:residue}.

\begin{prop} \label{prop:euler} Assume $T=S\times H$, $\pi: \lieT\to \lieS$ is the
natural projection, $p\in \mu(M^H)$, and $q=\pi(p)$. Then
$\kappa_p^H: H^*_T(M^H)\to H^*_H(M^H\red_p S)$ takes the
$T$-equivariant Euler class of the normal bundle of~$M^H$ onto the
$H$-equivariant Euler class of the normal bundle of $M^H\red_p S$
 in~$M\red_q S$.
\end{prop}

\begin{proof} The argument is almost identical to the one used
in~\cite{GK:residue}. Namely, let $Z=\mu^{-1}(p)\cap M^H$ and let
$\pr$ be the projection from $Z$ to $M^H\red_p S$ and $i:Z\into
M^H$. Then
$$
i^*\nu(M^H)=\pr^*\nu (M^H\red_p S).
$$
The proposition follows from functoriality of the Euler class as a
map from  oriented vector bundles to cohomology.
\end{proof}

\section{Residues}
\label{sec:residues}

In this section we define the residue operations and discuss their
basic properties.

\subsection{Residues of meromorphic 1-forms in one variable.}
Think of the Riemann sphere as the one point compactification
$\CC\cup\{\infty\}$ with the complex coordinate $z$ on $\CC$. Let
$f(z)$ be a meromorphic function on the Riemann sphere with values
in a topological vector space $V$ which can be written as the
finite sum
$$
f(z)=\sum_jg_j(z) e^{i\lambda_jz},
$$
where $g_j(z)$ are rational functions of $z$ and $\lambda_j\in
\RR$. Then in the case when all $\lambda_j\neq 0$ we define
$$
\res(f\ dz)=\sum_{\lambda_j>0}\sum_{b\in\CC} {\rm
res}(g_j(z)e^{i\lambda_jz};z=b)
$$
as  was done in~\cite{JK:LocalizationQuantization}. The other
case we will be interested in is when $\lambda_j=0$ for all $j$;
then we define
\begin{equation} \label{eq:resone}
\res (f\ dz)=\lim_{s\to 0^+} \res{f(z)e^{is\lambda z} dz}
\end{equation}
for some $\lambda>0$. It is easy to see that in this case $\res(f\
dz)$ is just the sum of all residues on $\CC$ and since the sum of
all residues of a meromorphic function is zero we conclude that
$$
\res(f\ dz)=-{\rm res}_{z=\infty}(f\ dz).
$$

Given a linear map $\psi:V\to W$ between two topological spaces,
the residue commutes with it
\begin{equation}
\label{eq:residue-commutes} \psi(\res(f)dz)=\res(\psi(f)dz).
\end{equation}

In the case when $V$ carries an algebra structure, an example
which will be important for us
is the residue of the function of the form
$$f=\frac{g(z)}{cz+a},
$$
where $a\in V$, $c\in \CC-\{0\}$, and $g(z)$ is a polynomial in
$z$ with values on $V$. For $A = \frac{a}{z}$ use
\begin{equation}
\label{eq:A} \frac{1}{1+A} = 1-A+A^2 - \dots,
\end{equation}
to  rewrite $f$ as a sum $\sum_{j = -
\infty}^{m_0} \gamma_j z^j $ for $\gamma_j \in V$, which  converges for
$|\frac{a}{z}|<1$.

Let $w=\frac{1}{z}$ be another coordinate on the Riemann sphere. Then
$$
\sum_{j = - m_0}^{\infty} \gamma_{-j} w^j  \left ( \frac{-dw}{w^2} \right )
$$
is the
Taylor
expansion of $fdz$ at $w=0$. Hence
\begin{equation}
\label{eq:-1} {\rm Res}(f\ dz)=-{\rm res}_{z=\infty}\Big(\sum_{j =
- \infty}^{m_0} \gamma_j z^j\ dz\Big) =
{\rm res}_{w=0}\Big(\sum_{j  =
- m_0}^{\infty} \gamma_{-j} w^{j-2}\ dw\Big) =
\gamma_{-1}.
\end{equation}

In the case when $g(z)$ is just a constant $g_0\in V$, we get
\begin{equation}
\label{eq:simple-residue} \res(f\ dz)=\frac{g_0}{c}.
\end{equation}

We emphasize that the residue is  well defined only as a function
of meromorphic 1-forms, not functions; the residue at $0 $ of the
1-form $ f(z) dz$ is independent of the choice of coordinate $z$
(invariant under a change of variables $z \mapsto g(z)$ provided
that $g(z)$ is a meromorphic function of $z$ and $dg(0) = 0 $
whereas this is not true of usual definition of the residue at $0$
of a meromorphic function).

\subsection{Residues of functions of several variables.}
Let us now consider a function $f$ of several complex variables
with values in a topological space $V$. More precisely, we assume
$f$ is  defined on the complexified Lie algebra
$\tc=\t\otimes\CC$ of the torus~$T$, and $f$ is a linear
combination of functions of the form
\begin{equation}
\label{eq:h} h(X)=\frac{q(X)e^{i\lambda(X)}}{\prod_{j=1}^k
\alpha_j(X)}
\end{equation}
for  some polynomials $q(X)$ of $X\in \t_\CC$, with values in $V$,
$\lambda\in \t^*$ and some $\alpha_1,\ldots,\alpha_k \in
\liets-\{0\}$.

Choose a coordinate system $X_1,\dots, X_m$ on $\t$, and denote by
the same letters the complexified coordinates which provide a
coordinate system of $\tc$. Let $\t_\ell$ be the subspace of $\t$
given by zeros of $X_1,\dots, X_{\ell}$. Define
$$
\res_{m}(h\ dX_m)=\res(h\ dX_m),
$$
where the variables $X_1,\ldots,X_{m-1}$ are held constant while
calculating this residue. As explained in Remark~3.5(1) of
\cite{JK:LocalizationQuantization}, in the case $\lambda\neq 0$
the residue $\res_m$ is well defined  only for a
generic  choice of coordinates $X_1,\dots, X_m$
(the precise condition on the coordinates
being specified in this Remark). Moreover, by
Remark~3.5(2) of \cite{JK:LocalizationQuantization}, $\res_{m}(f\
dX_m)$ is a linear combination of functions of the form
\begin{equation}
\label{eq:h2} \tilde h(X)=\frac{\tilde
q(X)e^{i\lambda(X)}}{\prod_{j=1}^{k-1} \tilde \alpha_j(X)}
\end{equation}
where $\tilde q(X)$ is a polynomial on the complexification of
$\t/\t_{m-1}$ and  $\tilde \alpha_j$ are in the dual of
$\t/\t_{m-1}$.

To consider the case when $\lambda = 0  $,
in other
words
\begin{equation}
\label{eq:hodef}
 h^0(X)=\frac{q(X)}{\prod_{j=1}^k
\alpha_j(X)},
\end{equation}
we may make a choice of  $\lambda^0 \in \t^*$ for
which $\lambda^0 (X) = \sum_{j=1}^m \lambda^0_j X_j $
with $\lambda^0_m > 0 $ and
define
${\res}_m (h^0 dX_m)$ as
\begin{equation} \label{eq:3.8B}{\res}_m (h^0 dX_m) =  \lim_{s \to 0^+}
{\res}_m (h^0 e^{is \lambda^0} dX_m). \end{equation}
We can easily check
that the residue ${\res_m} $ is a continuous function of $s \in
\RR^+$ and this limit exists, so we can define ${\res}_m (h
dX_m)$ even when $\lambda = 0 $. It is still true in this
situation that ${\res}_m (f dX_m)$ is a linear combination of
functions of the form  given in (\ref{eq:hodef}) with $\lambda=0$.

In the case when $V$ is an algebra, an important generalization of
(\ref{eq:simple-residue}) is the following.

\begin{lemma}
\label{lem:residue1}
$$
\res_m\Big(\frac{g_0}{\sum_i c_i X_i +v}\ dX_m\Big)
=\frac{g_0}{c_m},
$$
where $g_0,v\in V$ and $c_i\in\CC$.
\end{lemma}

\begin{proof} The proof follows from (\ref{eq:-1}) after setting
$c=c_m$ and $a=\sum_{i=1}^{m-1}c_i X_i +v$.
\end{proof}

\subsection{Iterated residues.} We again consider  functions ~$f$ which
are linear combinations of functions of  the form~(\ref{eq:h}). As
was just explained, $\res_m(fdX_m)$ is a linear combination of
functions of the form~(\ref{eq:h2}), which allows us to take the
residue of this function again. So we set $\res_m^m=\res_m$ and
by induction define the iterated residue
$$
\res_\ell^m(f\ [dX]^m_\ell) =\res(\res_{\ell+1}^m (f\
[dX]_{\ell+1}^m)\ dX_\ell),
$$
where $[dX]^m_\ell$ stands for the form $dX_\ell\wedge \dots
\wedge dX_m$ and, as above, the coordinates
$X_1,\ldots,X_{\ell-1}$ are held constant while calculating this
residue. Clearly $\res_\ell^ m(f\ [dX]^m_\ell)$ is a function on
the complexification of $\t/\t_{\ell-1}$.

In the case $V$ is an algebra, a generalization of
Lemma~\ref{lem:residue1} states the following.

\begin{lemma}
\label{lem:residue2} For a generic choice of coordinates
$X_1,\dots, X_m$
$$
\res_\ell^m\Big(
\frac{g_0}{\prod_{i=1}^{m-\ell+1}(\alpha_i(X)+v_i)}[dX]^m_\ell\Big)=
\frac{g_0}{{\rm det}(\tilde \alpha)}
$$
where $g_0,v_i\in V$, $\alpha_i\in \t^*-\{0\}$, and $\tilde
\alpha$ is the the matrix $\{a_{ij}\}_{1\leq i\leq m-\ell+1;
\ell\leq j\leq m} $, where $\alpha_i(X)=\sum a_{ij} X_j$.
\end{lemma}

\begin{proof}
The proof is analogous to the proof of property (iv) for the iterated
residue in  Proposition~3.4
of~\cite{JK:LocalizationQuantization}.
Notice that for
$h^0 $ defined in
(\ref{eq:hodef}) the definition of $\res_\ell^m (h^0
[dX]^m_\ell) $ is given in \cite{JK:LocalizationQuantization}
as the limit as $s \to 0^+$ of
 $\res_\ell^m (h^0e^{i s \lambda(X) }
[dX]^m_\ell) $ where $\lambda(X) = \sum_j \lambda_j X_j $ and
we choose $\lambda$ so that all the~$\lambda_j$  satisfy
$\lambda_j > 0$.
\end{proof}

Given a linear map $\psi:V\to W$ between two topological spaces,
the iterated residue commutes with this map:
\begin{equation}
\label{eq:residue-commutes3}
\psi(\res_\ell^m(f[dX]^m_\ell))=\res_\ell^m(\psi(f)[dX]^m_\ell).
\end{equation}

\subsection{Residues with respect to cones.} Iterated
residues depend on the choice of coordinates on~$\t$. Fix an inner
product on $\t$. Let us define residues which depend only on this
inner product and a choice of a certain cone $\dualcone$ in $\t$.

We introduce a function $f$ which is a linear combination of functions
of the form~(\ref{eq:h}).  We consider the set where none of the
functions $\alpha_j$ appearing in the denominators of functions
$h$ become zero, namely the set
\begin{equation} \label{eq:alcovs}
\{X\in\liet_\ell : \alpha_j(X) \neq 0, \text{ for all }\alpha_j
\}.
\end{equation}
Let $\Lambda$ be an open cone, which is a connected component of
this set.

Then for a generic choice of coordinate system
$X=(X_1,\ldots,X_m)$ on $\tc$ for which $(0,\ldots,0,1)\in
\Lambda$ define the residue with respect to the cone $\Lambda$ by
\begin{equation}
\label{eq:resdef} \res^{\Lambda}(h\ [dX])= \Jac \res_1^m(h\ [dX])
\end{equation}
where $[dX]=[dX]_1^m$ and $\Jac$ is the determinant of any
$(m)\times (m)$ matrix whose columns are the  coordinates of an
orthonormal basis of $\t$ defining the same orientation on $\t$ as
the chosen coordinate system.

To guarantee that $\res^\Lambda(h\ [dX])$ is well defined
(where $h$ is of the
form~(\ref{eq:h})), we need to make one additional
assumption: we assume that $\lambda$ is not in any proper subspace
of $\t$ spanned by some~$\alpha_i$'s. It was shown
in~\cite{JK:LocalizationQuantization} that under the above
assumptions~$\res^\Lambda(h\ [dX])$ is well defined, does not
depend on the choice of the coordinates but only on the choice of
the cone $\Lambda$ and the inner product on $\t$.

Originally, the residue $\res^\Lambda$ was introduced
in~\cite{JK:LocalizationNonabelian} as a generalization of a
certain integral over a vector space.
In~\cite[Proposition~3.4]{JK:LocalizationQuantization}, it was
shown that the definition~(\ref{eq:resdef}) coincides with the
original definition of $\res^\Lambda$
from~\cite{JK:LocalizationNonabelian}.
In~\cite[Proposition~3.2]{JK:LocalizationQuantization} it was
shown that certain properties together with linearity uniquely
define~$\res^\Lambda$. Let us recall these properties:
\begin{enumerate}
\item
\label{cond1} Let $\alpha_1,\dots,\alpha_v\in \Lambda^*$ be
vectors in the dual cone. Suppose that
 $\lambda$ is not in any cone of
dimension $m-1$ or less spanned by a subset of the $\{\alpha_i\}$.
 If $J=(j_1,\dots, j_m)$ is a multi-index and
$X^J=X_1^{j_1}\dots X_m^{j_m}$ then
$$
{\rm
Res}^\Lambda\left(\frac{X^Je^{i\lambda(X)}[dX]}{\prod_{i=1}^{v}\alpha_i(X)}\right)=0
$$
unless all of the
following properties are satisfied:
\begin{enumerate}
\item $\{\alpha_i\}_{i=1}^v$  span $\t^*$ as a vector space,
\item $v-(j_1+\dots +j_m)\geq m$,
\item $\lambda\in \langle \alpha_1,\dots, \alpha_v\rangle^+,$ the
positive span of the vectors $\{\alpha_i\}.$
\end{enumerate}

\item If properties (1)(a)-(c) above are satisfied, then
$$
{\rm Res}^\Lambda\left(\frac{X^Je^{i\lambda(X)}[dX]}
{\prod_{i=1}^{v}\alpha_i(X)}\right) = \sum_{k\geq 0}\lim_{s\to
0^+} {\rm Res}^\Lambda\left(\frac{X^J (i\lambda(X))^k
e^{is\lambda(X)}[dX]}{k!\prod_{i=1}^{v}\alpha_i(X)}\right),$$ and
all but one term in this sum are $0$ (the non-vanishing term being
that with $k = v-( j_1+\dots +j_m)     -m$).

\item \label{it:nonzeroresidue} The
residue is not identically $0$. If properties (1) $(a)-(c)$ are
satisfied with $\alpha_{1},\dots, \alpha_{m}$ linearly independent
in $\mathfrak{t}^*$, then

$$
{\rm Res}^\Lambda\left(\frac{e^{i\lambda(X)}[dX]}
{\prod_{i=1}^{m}\alpha_i(X)}\right) = \frac{1}{{\rm
det}(\overline{\alpha})},
$$
where $\overline{\alpha}$ is
 the nonsingular
matrix whose columns are the coordinates of $\alpha_{1},\dots,
\alpha_{m}$ with respect to any orthonormal basis of $\t$ defining
the same orientation.
\end{enumerate}

When $\lambda$ is of the form $\sum_{i=1}^k  s_{i} \alpha_{i}$
where fewer than $m$ of the $s_i$ are nonzero, then we define
\begin{equation}
\label{eq:prop4} {\rm Res}^\Lambda\left(\frac{e^{i\lambda(X)}[dX]}
{\prod_{i=1}^{v}\alpha_i(X)}\right) =\lim_{s\to 0^+} {\rm
Res}^\Lambda\left(\frac{e^{i(\lambda(X)+s\rho(X))}[dX]}
{\prod_{i=m}^{v}\alpha_i(X)}\right)
\end{equation}
where  $\rho \in \liet^*$ is chosen so that $\rho(\xi)>0$ for all
$\xi\in\Lambda$, and for small $s$, $\lambda + s \rho$ does not
lie in any cone of dimension $m-1$ or less spanned by a subset of
the $\{\beta_j\}.$

We will need another property of residues.

\begin{lemma} \label{lem:minus} Let $f(X_1,\dots, X_m)$ be a function on $\t$
given by a linear combination of  functions of the
form~(\ref{eq:h}) with $\lambda\neq 0$. Moreover, assume  that for every
set of values $a_1,\dots, a_{m-1}$ of the variables $X_1,\dots,
X_{m-1}$ the function $g(z)=f(a_1,\dots, a_{m-1}, z)$ is
holomorphic. Let $\Lambda$ be an appropriate choice of cone such
that $\res^\Lambda(f[dX])$ is defined, in particular $\Lambda$
contains the point $(0,\dots,0,1)$.  Then
\begin{equation}
\label{eq:opposite-cone}
\res^\Lambda(f[dX])=\res^{-\Lambda}(f[dY])
\end{equation}
where $Y_1,\dots, Y_m$ is a set of coordinates such that
$(0,\dots,0,1)\in -\Lambda$.
\end{lemma}

\begin{proof} For fixed $a_1,\dots, a_{m-1}$, let
$$
g(z)=\sum g_j(z) e^{i\lambda_j z}.
$$
Then
$$
\res(g(z) dz)={\rm res}^+(g(z) dz)=\sum_{\lambda_j>0}
\sum_{b\in\CC} {\rm res} (g_j(z) e^{i\lambda_jz};z=b).
$$

Define
$$
{\rm res}^-(g(z) dz)=\sum_{\lambda_j <0} \sum_{b\in\CC} {\rm res}
(g_j(z) e^{i\lambda_j z};z=b).
$$
Since $g(z)$ is holomorphic we have
\begin{equation}
\label{eq:sum=0}{\rm res}^+(g(z) dz)+{\rm res}^-(g(z) dz)=0.
\end{equation}

Since $\res^{-\Lambda}$ does not depend on the choice of
coordinates as long as $(0,\dots,0,1) \in-\Lambda$, we may choose
$Y_i=X_i$ for $1\leq i\leq m-1$, and $Y_m=-X_m$. Then
\begin{align*}
\res_m(f(Y_1,\dots, Y_m) dY_m)&=\res_m (-f(X_1,\dots, X_{m-1},
-X_m) dX_m)={\rm res}^- (-f(X_1,\dots, X_m)dX_m)\\
&=^\dagger {\rm res}^+ (f(X_1,\dots, X_m)dX_m) =
\res_m(f(X_1,\dots, X_m) dX_m),
\end{align*}
where ($\dagger$) holds because of~(\ref{eq:sum=0}).
Now~(\ref{eq:opposite-cone}) follows immediately from the definition
of $\res^\Lambda$ using iterated residues.
\end{proof}

\section{A generalization of Guillemin-Kalkman localization.}
\label{sec:GK}

In this section we discuss how to obtain the localization formula
of Guillemin and Kalkman~\cite{GK:residue} by applying ABBV
localization and then residue operations on symplectic cuts along
certain cones of dimension one. The same approach but with cones
of higher dimension and iterated residues provides a
generalization of Guillemin-Kalkman localization.

\subsection{Guillemin-Kalkman localization}  Let $(M, \omega)$ be
a compact symplectic manifold with an effective Hamiltonian~$T$
action. Let $\mu:M\to \t^*$ be the moment map. Pick a one
dimensional cone~$\cone$ transverse to~$\mu$. Assume the cone
$\cone$ is generated by a single weight $\beta$. If $\dim T=m$,
consider all $(m-1)$ dimensional walls of~$\mu$ which~$\cone$
intersects. Every such wall $W_i$ is an image $\mu(M_i)$ of a
connected component $M_i$ of $M^{H_i}$ for some one dimensional
subtorus $H_i$ of $T$. If $p_i$ is the intersection $W_i\cap
\cone$, then $G_i=M_i\red_{p_i} T$ is a connected component of the fixed
point set of the $T/T^\cone$ action on~$M_\cone$. The only other
connected component of this fixed point set is $M_0$, the
reduction of $M$ at zero by~$T$. (If $\dim T=1$, then $G_i$ are
also connected components of the  fixed point set $M^T$, and
using the notation introduced in
Definition~\ref{def:points}, $G_i$ would be called
the old connected components of the fixed point set. If $\dim
T>1$, then in terms of  the same notation $G_i$ are the new fixed points.)

Recall that any equivariant cohomology class $\eta\in H^*_T(M)$
descends to a class $\eta_\cone$ on $M_\cone$:
$$
\eta_\cone=\kappa_\cone(\eta\otimes 1)\in
H^*_{T/T^\cone}(M_\cone)\otimes 1\subset H^*_T(M_\cone).
$$
Think of $\eta_\cone$ as a class in $H^*_{T/T^\cone}(M_\cone)$ and
apply the orbifold version of the ABBV localization theorem to it:
\begin{equation}
\label{eq:ABBV-on-cut}
\frac{1}{d_{M_\cone}}\int_{M_\cone}\eta_\cone=
\frac{1}{d_{M_0}}\int_{M_0}\frac{\kappa_0(\eta)}{e(\nu(M_0))}
+\sum_{G_i} \frac{1}{d_{G_i}}\int_{G_i}
\frac{\iota_{G_i}^*(\eta_\cone)}{e(\nu(G_i))}.
\end{equation}
where, as usual, $\kappa_0$ is the Kirwan map at zero and
$\iota_{G_i}$ are the inclusions of $G_i$ into $M^\cone$.

Apply the residue operation to both sides
of~(\ref{eq:ABBV-on-cut}). For this we have to specify  a
coordinate $X_1$ on the one dimensional Lie algebra of
$T/T^\cone$. Clearly the dual of this Lie algebra can be
identified with the line in $\t^*$ passing through the cone
$\cone$. So $\beta$ defines a coordinate on the Lie algebra of
$T/T_\cone$, and we define $X_\beta=X_1$ to be this coordinate.

Using $X_1$ we get the formula
\begin{equation}
\label{eq:GK1}
\res_1\frac{dX_\beta}{d_{M_\cone}}\int_{M_\cone}\eta_\cone=
\res_1\frac{dX_\beta}{d_{M_0}}\int_{M_0}\frac{\kappa_0(\eta)}{e(\nu(M_0))}
+\sum_{G_i} \res_1\frac{dX_\beta}{d_{G_i}}\int_{G_i}
\frac{\iota_{G_i}^*(\eta_\cone)}{e(\nu(G_i))}.
\end{equation}
Since $\int_{M_\cone}\eta_\cone$ is just a polynomial function on
the Lie algebra of $T/T^\cone$, by~(\ref{eq:-1}) the left hand
side of (\ref{eq:GK1})~is zero.

The normal bundle $\nu(M_0)$ is just a line bundle whose Euler
class is given by
\begin{equation}
\label{eq:Euler}
e(\nu(M_0))=X_\beta +c_1(\nu(M_0)).
\end{equation}
where $c_1(\nu(M_0))$ is the first Chern class of $\nu(M_0)$. Thus
by~(\ref{eq:residue-commutes3}) and~(\ref{eq:simple-residue})
$$
\res_1\frac{dX_\beta}{d_{M_0}}\int_{M_0}\frac{\kappa_0(\eta)}{e(\nu(M_0))}
=^\dagger\frac{1}{d_{M_0}}\int_{M_0}\res_1\Big(\frac{\kappa_0(\eta)\
dX_\beta }{X_\beta +c_1(\nu(M_0))}\Big) =\frac{1}{d_{M_0}}
\int_{M_0}\kappa_0(\eta).
$$

\begin{remark}
\label{rem:commute}
To show that ($\dagger$) follows from~(\ref{eq:residue-commutes3}) think
of  $\frac{\kappa_0(\eta)}{e(\nu(M_0))}$ as a function of $\t$ with values
in $H^*(M)$. If we define $\psi:H^*(M)\to\CC$ to be the usual
integration on $M$, then~(\ref{eq:residue-commutes3}) applied to this
$\psi$ yields~($\dagger$).
\end{remark}

Hence (\ref{eq:GK1}) yields a formula for the the integral of
$\kappa_0(\eta)$ over $M_0$ of the Kirwan map of $\eta$ in terms
of the residues of certain integrals over the $G_i$. To put this
formula in the form in which it appeared in~\cite{GK:residue} let us
transform the terms corresponding to the $G_i$ in~(\ref{eq:GK1})
by commuting $\res_1$ with the integration and the Kirwan map.

By~(\ref{eq:residue-commutes3}) and Proposition~\ref{prop:euler} we have
$$
\res_1\Big(\frac{dX_\beta}{d_{G_i}}\int_{G_i}
\frac{\iota_{G_i}^*(\eta_\cone)}{e(\nu(G_i))}\Big)=
\frac{1}{d_{G_i}}\int_{G_i}\res_1 \Big(\kappa_{p_i}^{H_i}
\Big(\frac{\iota^*_{M_i}\eta}{e(\nu(M_i))}\Big)dX_\beta\Big)
$$
where $\iota_{M_i}$ is just the inclusion of $M_i$ into $M$, and
as defined in~(\ref{eq:kappa_H}) $\kappa_{p_i}^{H_i}$ is the
Kirwan map from $H^*_T(M_i)$ to $H^*_{T/T^\cone}(G_i)$.

Notice that each $M_i$ is fixed by the one-dimensional torus
$H_i$, so that the Lie algebra~$\mathfrak h_i$ is a subspace of
$\t$. Pick a coordinate system $X^i_1,\dots,X^i_m$ on $\t$ such
that $X^i_m=\beta$ and $X^i_1,\dots, X^i_{m-1}$ vanish
on~$\mathfrak h_i$, so that $\t_{m-1}=\mathfrak h_i$. Then it is
not difficult to see that
$$
\res_1 \circ \kappa_{p_i}^{H_i} =\kappa_{p_i}^{H_i} \circ
\res_{m,i},
$$
where $\res_{m,i}$ is $\res_m$ defined using the coordinate
system $\{X_j^i\}$ (note that the definition of
$\res_m$ was given in
(\ref{eq:3.8B}).
 This allows us to obtain the following restatement of
\cite[Theorem~3.1]{GK:residue}.

\begin{theorem}
\label{thm:GK} For $\eta\in H^*_T(M)$
\begin{equation}
\label{eq:thmGK}
\int_{M_0}\kappa_0(\eta)=-\sum_{G_i}\frac{d_{M_0}}{d_{G_i}}
\int_{G_i}\kappa_{p_i}^{H_i} \res_{m,i}
\Big(\frac{\iota^*_{M_i}\eta}{e(\nu(M_i))}d\beta\Big).
\end{equation}
\end{theorem}

Guillemin and Kalkman iterated this result using certain
combinatorial objects called dendrites to produce a formula for
$\int_{M_0}\kappa_0(\eta)$ in terms of the integration over the
connected components of the fixed point set $M^T$ of the original
$T$ action on $M$. We will discuss dendrites and their
generalizations in Section~\ref{sec:dendrites}.

\subsection{Higher dimensional generalization of Guillemin-Kalkman
localization.} We now generalize the formula of the
previous section to the case when we cut with a cone of any
dimension.

As before, let $(M, \omega)$ be a compact symplectic manifold with
an effective Hamiltonian $T$ action. Let $\mu:M\to \t^*$ be the
moment map. Pick any cone $\cone$ transverse to $\mu$. Assume
$\cone$ is generated by  the linearly independent weights
$\beta_1,\dots, \beta_k$. Recall that in
Definition~\ref{def:points} we distinguished three kinds of fixed
points of $M_\cone$: the old fixed points (whose connected
components are denoted by $F_i$), the new fixed points (whose
connected components are denoted by $F_i'$), and the symplectic
reduction~$M_0$.

%Let us consider the intersections of the walls of $\mu(M)$ and the
%open faces of $\cone$ of complementary dimensions, so that they
%intersect along points. These points correspond to the connected
%components of the fixed point set of the $T$ action on $M_\cone$.
%Let us distinguish three kinds of such points: the points $p_i$,
%which are the intersections of the zero dimensional walls of
%$\mu(M)$ and the interior of $\cone$ (these points exist only in
%the case $\dim\cone =\dim \lieT$); the points $p_i'$, which are
%the intersections of a wall of $\mu(M)$ of dimension  $r
%> 0 $ and a $\dim\lieT -r$ dimensional open face of $\cone$; and
%the point $0$. These three sets of points are the images under the
%map $\mu_{\rm cut}$ of the corresponding components of
%$M^T_\cone$, as defined in Section~\ref{subsec:cut}.

As in the case of one-dimensional cones any equivariant cohomology class
$\eta\in H^*_T(M)$ descends to a class $\eta_\cone$ on $M_\cone$
$$
\eta_\cone=\kappa_\cone(\eta\otimes 1)\in
H^*_{T/T^\cone}(M_\cone)\otimes 1\subset H^*_T(M_\cone).
$$
As before,  we think of $\eta_\cone$ as a class in
$H^*_{T/T^\cone}(M_\cone)$ and apply the ABBV localization theorem to
it:
\begin{equation}
\label{eq:ABBV-cones}
\frac{1}{d_{M_\cone}}\int_{M_\cone}\eta_\cone=
\frac{1}{d_{M_0}}\int_{M_0}\frac{\kappa_0(\eta)}{e(\nu(M_0))}
+\sum_{F_i} \frac{1}{d_{F_i}}\int_{F_i}
\frac{\iota_{F_i}^*(\eta_\cone)}{e(\nu(F_i))} +\sum_{F'_i}
\frac{1}{d_{F'_i}}\int_{F'_i}
\frac{\iota_{F'_i}^*(\eta_\cone)}{e(\nu(F'_i))}.
\end{equation}
Here,  $\kappa_0$ is the Kirwan map at zero and
$\iota_{F_i},\iota_{F_i'}$ are the inclusions of $F_i, F_i'$ into
$M_\cone$.

Let us apply the iterated residue operation to both sides
of~(\ref{eq:ABBV-cones}). As in the one-dimensional case the dual
of the Lie algebra of $T/T^\cone$ can be identified with the
subspace in $\t^*$ passing through the cone~$\cone$. Consider the
set of weights $\{\alpha_i\}$ which appear in the isotropy
representation of the $T/T^\cone$ action at the connected
components of the fixed point set of the $T/T_\cone$ action on
$M_\cone$. Notice that the weights $\tilde \beta_1,\dots,\tilde
\beta_k$ of the isotropy representation of $T/T^\cone$ at $M_0$
are the images of the
weights $\beta_1,\dots,\beta_k$
under the map $(\t/\t^\cone)^*\to\t^*$. In particular, the
set
$\{\alpha_i\}$ contains the weights $\tilde \beta_i$. Let
$\{X_1,\dots,X_k\}$ be a generic coordinate system on
$T/T^\cone$ such that the determinant of the matrix %with coefficients
which expresses $\tilde \beta_i$ as a linear combination of $X_j$ is one.

%Pick
%$\Lambda$ to be any connected components of the set
%$$
%\{X\in{\rm Lie}(T/T_\cone)| \alpha_i(X)\neq 0\}
%$$
%containing points which evaluate nonnegatively on the points of
%the cone $\cone$. In other words, the dual cone $\cone^*$ must
%contain $\Lambda$.

Then  by applying $\res_1^k$ to both sides
of~(\ref{eq:ABBV-cones}) we get the formula
\begin{equation}
\label{eq:GKmulti} 0= \res_1^k\frac{[dX]_1^k}{d_{M_0}}
\int_{M_0}\frac{\kappa_0(\eta)}{e(\nu(M_0))} +\sum_{F_i}
\res_1^k\frac{[dX]^k_1}{d_{F_i}}\int_{F_i}
\frac{\iota_{F_i}^*(\eta_\cone)}{e(\nu(F_i))} + \sum_{F'_i}
\res_1^k\frac{[dX]^k_1}{d_{F'_i}}\int_{F'_i}
\frac{\iota_{F'_i}^*(\eta_\cone)}{e(\nu(F'_i))},
\end{equation}
where the left hand side is zero, since $\int_{M_\cone}\eta_\cone$
is just a
polynomial function on the Lie algebra of~$T/T^\cone$.

The weights of the isotropy representation of $T/T_\cone$ on the
normal bundle $\nu(M_0)$ are just $\tilde \beta_i$. Hence
$\nu(M_0)$ splits as a direct sum of line bundles $\oplus_i L_i$,
where $T/T^\cone$ acts on the fibers of $L_i$ by the weights
$\tilde \beta_i$. Then the Euler class of $\nu(M_0)$ is just
\begin{equation}
\label{eq:Euler2} e(\nu(M_0))=\prod_i(\tilde\beta_i +c_1(L_i)).
\end{equation}
where $c_1(L_i)$ are the first Chern classes of $L_i$. Thus
by~(\ref{eq:residue-commutes3}) (see also Remark~\ref{rem:commute}) and
Lemma~\ref{lem:residue2} we have
\begin{equation} \label{eq:foursixb}
\res_1^k\frac{[dX]_1^k}{d_{M_0}}
\int_{M_0}\frac{\kappa_0(\eta)}{e(\nu(M_0))} =\frac{1}{d_{M_0}}
\int_{M_0}\kappa_0(\eta).
\end{equation}

By analogy with the one-dimensional case, we should now simplify
the second and third terms of the right hand side
of~(\ref{eq:GKmulti}). But the second term is already in the form
we want, since it is written in terms of the fixed points
of $M$.

To simplify the last term of~(\ref{eq:GKmulti}), first consider
the situation when a point $p_i'=\mu_\cone(F_i')$ is an intersection of
the
interior of $\cone$ with a wall of $\mu(M)$ of dimension $r>0$.
Assume that this wall is the moment map image of $M_i$, which is
stabilized by $H_i\subset T$. Then by~(\ref{eq:residue-commutes3}) and
Proposition~\ref{prop:euler}
we have
\begin{equation}
\label{eq:new-fixed}
\frac{1}{d_{F'_i}}\res_1^k\Big([dX]_1^k\int_{F'_i}
\frac{\iota_{F'_i}^*(\eta_\cone)}{e(\nu(F'_i))}\Big)=
\frac{1}{d_{F'_i}}\int_{F_i}\res_1^k
\Big([dX]_1^k\kappa_{p'_i}^{H_i}
\big(\frac{\iota^*_{M_i}\eta}{e(\nu(M_i))}\big)\Big)
\end{equation}
where $\iota_{M_i}$ is just the inclusion of $M_i$ into $M$, and
as defined before $\kappa_{p'_i}^{H_i}$ is the Kirwan map from
$H^*_T(M_i)$ to $H^*_{T/T^\cone}(F'_i)$.

For the more general situation when $p_i$ is an intersection of
any face of $\cone$ with a wall of $\mu(M)$ of complementary
dimension, formula~(\ref{eq:new-fixed}) becomes a little more
complicated. Namely the Euler class of~$F_i'$ is no longer a
Kirwan map image of the Euler class of $M_i$, but rather this
multiplied by another class $\nu_i$. The class $\nu_i$ can be
understood as follows. Take the symplectic cut $M_{\cone,i}$
of~$M_i$ with respect to the cone $\cone$. This symplectic cut
sits naturally  inside $M_\cone$. Then $\nu_i$ is just the Euler
class of the normal bundle of $F_i'$ inside $M_{\cone,i}$. While
this class can be quite complicated, by the Kirwan surjectivity
theorem we know that there exists a class $\tau_i$ on $M_i$ such
that $\kappa_{p_i'}^{H_i}(\tau_i)=\nu_i$. So we conclude that
$$
\frac{1}{d_{F'_i}}\res_1^k\Big([dX]_1^k\int_{F'_i}
\frac{\iota_{F'_i}^*(\eta_\cone)}{e(\nu(F'_i))}\Big)=
\frac{1}{d_{F'_i}}\int_{F'_i}\res_1^k\Big([dX]_1^k
\kappa_{p'_i}^{H_i}
\big(\frac{\iota^*_{M_i}\eta}{e(\nu(M_i))\tau_i}\big)\Big).
$$

To simplify this term even further, let us choose a  set of
coordinates $\{ X^i_1,\dots, X_m^i\}$ on $\t$ as follows. Each~$M_i$ is
fixed by a subtorus $H_i$ of $T$, so that its Lie algebra
$\mathfrak h_i$ lies inside $\t$. (Notice that $\dim H_i\leq \dim
\cone$ and this inequality might be strict.). Choose any generic
coordinates $X_1^i,\dots X^i_{m-k}$ which vanish on $\mathfrak
h_i$ and coordinates $X^i_{m-k+1}, \dots, X^i_m$, which are the images
of $X_1, \dots, X_k$ under the map $(\t/\t^\cone)^*\to \t^*$.

Then we
conclude that
$$
\res_1^k \circ \kappa_{p'_i}^{H_i} =\kappa_{p'_i}^{H_i} \circ
\res_{m-k+1}^m,
$$
where $\res_{m-k+1}^m$ on the right hand side is take with respect to
coordinates $X_1^i,\dots,X_m^i$.
Thus we get the following generalization of
\cite[Theorem~3.1]{GK:residue} to the case of cones of arbitrary
dimension.

\begin{theorem}
\label{thm:GK-multi} For $\eta\in H^*_T(M)$
\begin{equation}
\label{eq:thmGK-multi} \int_{M_0}\kappa_0(\eta)=-\sum_{F_i}
\frac{d_{M_0}}{d_{F_i}}\int_{F_i} \res_1^k[dX]_1^k
\frac{\iota_{F_i}^*(\eta_\cone)} {e(\nu(F_i))}-
\sum_{F'_i}\frac{d_{M_0}}{d_{F'_i}} \int_{F'_i}\kappa_{p'_i}^{H_i}
\res_{m-k+1}^m
[dX^i]^m_{m-k+1}\frac{\iota^*_{M_i}\eta}{e(\nu(M_i))\tau_i}.
\end{equation}
\end{theorem}

\subsection{Dendrites and their generalizations}
\label{sec:dendrites}

In addition to  Theorem~\ref{thm:GK}, Guillemin and
Kalkman~\cite{GK:residue} proved a formula which uses
Theorem~\ref{thm:GK}
 iteratively to express $\int_{M_0}\kappa_0(\eta)$ in terms of
integration over the connected components of $M^T$. They called
the combinatorial objects responsible for how the iteration is
performed dendrites.

Consider a finite collection $D$ of tuples $(\cone(q), W)$, where
$\cone(q)$ is a one-dimensional shifted cone (that is, if
$\cone\subset \t^*$ is a cone centered at the origin and
$q\in\t^*$, then $\cone(q)= \{x+q\ |x\in \cone\}$) and $W$ is a
wall of $\mu(M)$ such that $\cone(q)$ lies inside the affine
space spanned by $W$. Notice that according to our conventions
every $W$ uniquely defines a subtorus $H_W$ of $T$ and a connected
component $M_W$ of $M^{H_W}$ such that $W=\mu(M_W)$. Assume that
$D$ contains a tuple $(\cone_0(0), \mu(M))$ (we will treat
$\cone_0$ as the cone used in Theorem~\ref{thm:GK}). Then we say
that  $D$ is a dendrite if the following conditions hold
\begin{enumerate}
\item For $(\cone(q), W)\in D$ the cone $\cone(q)$ is transverse to
$\mu(M_W)$,
\item For $(\cone(q), W)\in D$ if $\cone(q)$ intersects a codimension one
wall $W'$ of
$\mu(M_W)$ at a point $q'$, then there is a unique cone $\cone'$ such that
$(\cone'(q'),W')\in D$.
\end{enumerate}
Obviously it is possible to construct a dendrite which contains a tuple
$(\cone_0(0), \mu(M))$ as long as $\cone_0$ is transverse
to~$\mu(M)$.

Now consider formula~(\ref{eq:thmGK}) of Theorem~\ref{thm:GK}
produced using the cone $\cone_0$. Assume we are given a dendrite
$D$ containing $(\cone_0(0),\mu(M))$. Every term in the summation
on the right hand side of~(\ref{eq:thmGK}) corresponds to an
intersection $p_i=\mu_\cone(G_i)$ of $\cone_0$ with a codimension one wall
$W_i$.
By definition of the dendrite, there exists a unique cone
$\cone_i$ such that $(\cone_i(p_i),W_i)\in D$.  Set $H_i=H_{W_i}$
and $M_i=M_{W_i}$ and notice that  $G_i=M_i\red_{p_i}
(T/H_i)$. Moreover, let us shift the moment map on $M_i$ by
defining  $\mu_i=\mu-p_i$. Then, we can apply Theorem~\ref{thm:GK}
to the $T/H_i$ action on $M_i$ and the cone $\cone_i$  to compute
$$
\int_{M_i\red_{p_i} (T/H_i)}\kappa_{p_i}^{H_i} \res_m \Big(dX_m
\frac{\iota^*_{M_i}\eta}{e(\nu(M_i))}\Big)
$$
as a sum of integrals over  symplectic reductions of the form
$N\red (T/H)$, where $H$ is a two dimensional subtorus of $T$ and
$N$ is a connected component of $M^H$.
%Indeed, $\res_m \Big([dX_m]
%\frac{\iota^*_{M_i}\eta}{e(\nu(M_i))}\Big)$ is a $T/H_i$
%equivariant cohomology class on $M_i$.

Repeating this process we can express $\int_{M_0}\kappa_0(\eta)$
as a summation of integrals over the connected components of
the fixed point set $M^T$.  More specifically, for a dendrite $D$
we say that the sequence of tuples from $D$
$$\big( (\cone_{0}(q_{0}),W_{0}),
(\cone_{1}(q_{1}),W_{1}), \dots, (\cone_{k}(q_{k}),W_{k})\big)$$
is a path $P$ if $(\cone_0,\mu(M))=(\cone_{0}(q_{0}),W_{0})$;
$W_k$ is zero-dimensional, so that there is a connected component
$F$ of $M^T$ with $\mu(F)=W_k$ (we will denote this $F$ by $F_P$);
each $q_{j}$ is an intersection of $\cone_{j-1}(q_{j-1})$ with a
codimension one wall $W_{j}$ of $\mu(M_{W_{j-1}})$. For a path
$P$, let $k_{j}=\dim W_{j}$. Then for $\gamma\in
H_{T/H_{W_{j-1}}}^*(M_{W_{j-1}})$ and an appropriately chosen set
of coordinates on the Lie algebra of $T/H_{W_{j-1}}$ define
$$
Q_j (\gamma)=-\res_{k_{j}}\Big(dX_{k_j}
\frac{\iota^*_{M_{W_j}}\gamma}{e(\nu(M_{W_j}))}\Big).
$$
For a path $P$ define
$$
Q_P=Q_k\circ\dots\circ Q_1.
$$
The above discussion proves the following.

\begin{theorem}
\label{thm:dendrite} For a dendrite $D$ and $\eta\in H^*_T(M)$
$$
\int_{M_0}\kappa_0(\eta)=d_{M_0}\sum_P \int_{F_P} Q_P(\eta)
$$
where the above sum is taken over all possible paths in~$D$.
\end{theorem}

%Moreover, as shown
%in~\cite{GK:residue} it is possible to significantly simplify each
%term of the summation corresponding to $F_i$, and write it as an
%integral of a series of residues applied to
%$\frac{\iota_{F_i}^*\eta}{e(\nu(F_i))}$. We refer the reader
%to~\cite{GK:residue} for details and do not reproduce them, since
%as will be explained later these simplifications of terms do not
%generalize to the multi-dimensional dendrites.

Now consider a collection $D$ of tuples $(\cone(q), W)$, where
$\cone(q)$ is a shifted cone of any dimension and~$W$ is a wall of
$\mu(M)$ such that $\cone(q)$ lies inside the affine space spanned
by $W$. Again, every~$W$ uniquely determines a subtorus $H_W$ of
$T$ and a connected component $M_W$ of $M^{H_W}$ such that
$W=\mu(M_W)$. Assume that $D$ contains a tuple $(\cone_0(0),
\mu(M))$. Then we say that  $D$ is a multi-dimensional dendrite if
the following conditions hold
\begin{enumerate}
\item For $(\cone(q), W)\in D$ the cone $\cone(q)$ is transverse to
$\mu(M_W)$,
\item For $(\cone(q), W)\in D$ if a face of $\cone(q)$ of dimension $r\leq \dim W$
intersects a wall $W'$ of $\mu(M_W)$ of dimension $\dim W-r$ at a
point $q'$, then there is a unique cone $\cone'$ such that
$(\cone'(q'),W')\in D$.
\end{enumerate}
%Again, it is possible to construct a dendrite which contains a
%tuple $(\cone_0(0), \mu(M))$ as long as $\cone_0$ is transverse
%to~$\mu(M)$.

In analogy with  the case of dendrites of one-dimensional cones,
consider formula~(\ref{eq:thmGK-multi}) of
Theorem~\ref{thm:GK-multi} produced using the cone $\cone_0$.
Assume we are given a multi-dimensional dendrite $D$ containing
$(\cone_0(0), \mu(M))$. Every term in the second summation on the
right hand side of~(\ref{eq:thmGK-multi}) corresponds to a new
connected component $F'_i$ of $M_\cone^T$. Moreover, its moment
map image $p'_i=\mu_\cone(F'_i)$ is always an intersection of a
face of $\cone$ with a wall $W_i$ of complementary dimension. By
definition of the multi-dimensional dendrite, there exists a
unique cone $\cone_i$ such that $(\cone_i(p_i'),W_i)\in D$.
(Notice that $M_{W_i}$ and $H_{W_i}$ are just $M_i$ and $H_i$ in
the notation of Theorem~\ref{thm:GK-multi}.) As in the case of
one-dimensional cones, we know that $F'_i=M_{i}\red_{p_i}
(T/H_{i})$, so we can apply Theorem~\ref{thm:GK-multi} to the
$T/H_{i}$ action on $W_i$ and the cone $\cone_i(p'_i)$  to compute
$$
\int_{M_{i}\red_{p_i} (T/H_{i})}\kappa_{p'_i}^{H_{i}}
\res_{m-k+1}^m
\Big([dX^i]^m_{m-k+1}\frac{\iota^*_{M_i}\eta}{e(\nu(M_i))\tau_i}\Big)
$$
as a sum of integrals over  symplectic reductions of the form
$N\red (T/H)$, where $H$ is a subtorus of
$T$ of dimension at least $2$ and $N$ is a connected component of $M^H$.

Using the multi-dimensional dendrite we can iterate the process
and express $\int_{M_0}\kappa_0(\eta)$ as a summation of integrals
over the connected components of $M^T$ as follows. For a
multi-dimensional dendrite $D$ we say that the sequence of tuples
$$\big( (\cone_{0}(q_{0}),W_{0}),
(\cone_{1}(q_{1}),W_{1}), \dots, (\cone_{k}(q_{k}),W_{k})\big)$$
is a path $P$ if $(\cone_0,\mu(M))=(\cone_{0}(q_{0}),W_{0})$;
$W_k$ is zero-dimensional, so that there exists a connected
component $F$ of $M^T$ with $\mu(F)=W_k$ (again we denote
$F_P=F$); each $q_{j}$ is the intersection of the interior of a
face $\sigma$ of $\cone_{j-1}(q_{j-1})$ with a wall $W_{j}$ of
$\mu(M_{W_{j-1}})$ such that $\dim\sigma+\dim W_{j} =\dim
W_{j-1}$. For a path $P$, let $k_{j}=\dim W_{j}$  and $m_j=\dim
\cone_j$. Then for $\gamma\in H_{T/H_{W_{j-1}}}^*(M_{W_{j-1}})$
and an appropriately chosen set of coordinates on the Lie algebra
of $T/H_{W_{j-1}}$ define
$$
Q_j (\gamma)=-\res_{k_j-m_j+1}^{k_j}\Big([dX]^{k_j}_{k_j-m_j+1}
\frac{\iota^*_{M_{W_j}}\gamma}{e(\nu(M_{W_j}))}\Big).
$$
For a path $P$ define
$$
Q_P=Q_k\circ\dots\circ Q_1.
$$
It is clear that Theorem~\ref{thm:dendrite} holds in the case when
$D$ is a multi-dimensional dendrite.

\section{A new proof of the Jeffrey-Kirwan localization formula.}
\label{sec:JK} As explained in the previous section it is possible
to iterate the formulas of Theorems~\ref{thm:GK}
and~\ref{thm:GK-multi} using (multi-dimensional) dendrites to
express integration of $\kappa_0(\eta)$ over $M_0$ as a sum of
integrals of certain forms over the connected components of $M^T$.
In this section we present another way of writing such a formula.
The rough idea is to choose a very wide cone $\cone$, whose dual
is inside a certain cone $\Lambda$, then apply the residue
operation with respect to $\Lambda$  to the ABBV formula for $\eta
e^{i\tilde \omega}$ (where $\tilde \omega$ is the equivariant
symplectic form) in such a way that  the terms corresponding to
the new fixed points $F_i'$ of $M_\cone$ are zero. As a result we
will obtain a new proof of the the Jeffrey-Kirwan localization
theorem~\cite{JK:LocalizationNonabelian}.

\subsection{The choice of the cone}  In this section we explain how to
choose the cone $\cone$. It is convenient to think of $\cone$ as a
very wide cone, in other words a cone which is close to being a hyperplane.

Given an effective Hamiltonian action
of $T$ on a compact manifold~$M$, let $\{\alpha_i\}$ be the set of
all weights appearing in the isotropy representation of $T$ at the
fixed points. Pick a connected component $\Lambda$ of the set
$$
\{\xi\in \t | \alpha_i(\xi)\neq 0\text{ for all }i\}.
$$
Consider the dual cone  $\Lambda^*=\{X\in\t^*| X(\xi)\geq
0\text{ for all }\xi \in \Lambda\}$.
%\comment{Misha -- I removed ``dual cone of'' above -- Lisa}
Now pick a cone $\cone$ transverse to~$\mu(M)$ spanned by $m=\dim
\t$ weights $\beta_1,\dots,\beta_m$  such that $\cone$ satisfies
\begin{enumerate}
\item
$\Lambda^*\subseteq \cone$, or in other words $\cone^*\subseteq
\Lambda$,
\item
For every wall $W$ of $\mu(M)$ if $W\cap \cone$ is not empty,
then there exists
$\xi_W\in \cone^*$ such that the maximum of
$q(\xi_W)$ for  $q\in W\cap \cone$ is
attained at a vertex of $W$.
\end{enumerate}

Let us explain why such a cone $\cone$ always exists. Pick a
rational vector $\xi_0\in\Lambda$ and let $H$ be the hyperplane in
$\t^*$ annihilated by $\xi_0$. For a weight $p$  in the interior
of $ \Lambda^*$ define $H_p= H+p$. Then $H_p\cap \Lambda^*$ is a
rational polytope. In particular, there exists a simplex $S\subset
H_p$ with rational vertices $\tilde \beta_1,\dots, \tilde \beta_m$
which contains the polytope $H_p\cap \Lambda^*$. Take the weights
$\beta_1,\dots, \beta_m$, which define $\cone$, to be multiples of
$\tilde \beta_1,\dots, \tilde \beta_m$. Clearly $\cone$ satisfies
(1). Moreover, by increasing the size of $S$, that is by making
$\cone$ wider, we can always guarantee that $\cone $ satisfies
(2). Since there are infinitely many choices of simplices $S$,
there is a choice of $S$ for which $\cone$ is transverse to
$\mu(M)$.

Now consider the symplectic cut $M_\cone$ and the set
$\{\gamma_i\}$ of all  the weights  which appear as weights of
isotropy representations at  fixed points of $M_\cone$.
%(Notice that the weights $\beta_i$ are in the set $\{\gamma_i\}$.)
Consider the set
$$
\cone^*\cap\{\xi\in \t | \gamma_i(\xi)\neq 0\text{ for all }i\},
$$
and let the cone $\bar \Lambda$  be a connected component of this
set.

For each connected component $F$ of $M_\cone^T$, let
$\{\gamma_j^F\}$ be the set of weights of the isotropy
representation of $T$ at $F$. Polarize the weights
$\{\gamma_j^F\}$ using $\bar \Lambda$, that is if
$\gamma^F_j(\xi)<0$ for all $\xi\in\bar \Lambda$ then set
$\bar\gamma^F_j=\gamma^F_j$, otherwise let
$\bar\gamma^F_j=-\gamma^F_j$. Define $C_F$ to be the cone
containing all the points of the form $\mu_\cone(F)+\sum s_j
\bar\gamma_j^F$ with $s_j$ being nonnegative real numbers.

If $F$ is an old connected component of the fixed point set, then
the polarization of the weights $\{\gamma_j^F\}$ using $\bar
\Lambda$ is the same as using $\Lambda$, since $\bar
\Lambda\subset \Lambda$. If $F=M_0$, then $C_F=-\cone$. For the
new fixed points let us prove the following important fact.

\begin{lemma}
\label{lem:polarize} Let $F$ be a new connected component of the
fixed point set $M_\cone^T$. Then the cone $C_F$ does not
intersect the interior of $\cone$.
\end{lemma}

\begin{proof} Let $p=\mu_{\cone}(F)$. Since $\cone$ intersects $\mu(M)$
transversely, there is a subtorus $H$ of $T$ and a connected
component $M'$ of $M^H$ for which  $F$ is the symplectic reduction
of $M'$ at $p$. Moreover,~$p$ is just the intersection of the wall
$W=\mu(M')$ of dimension $k$ and an open face $\sigma$ of $\cone$
of dimension~\mbox{$m-k$.}

Let us order the weights $\{\gamma_j^F\}_{j=1}^n$ of the isotropy
representation of $T$ at $F$ in a special way. Assume that
each of the first $k'$ weights is parallel to an
intersection $W\cap \sigma'$, where $\sigma'$ is an open face of
$\cone$ such that $\dim \sigma'=\dim \sigma+1$ and the closure of
$\sigma'$ contains $\sigma$. Assume  also that each of the last
$n-k'$ weights is parallel to an intersection
$W'\cap \sigma$, where $W'$ is a wall with $\dim W'=\dim W+1$ and
$W\subset W'$. (It is easy to see that every $\gamma_j^F$ falls
into one of these  categories.)

Let $C_F^1$ be the cone spanned by the first $k'$ polarized
weights centered at the origin and let $C_F^2$ be the cone spanned by
the other weights, also centered at the origin. Then
$C_F=p+C_F^1+C_F^2$.  Since the cone $p+C_F^2$ lies in the affine
space passing through $\sigma$, it remains to show that $p+C_F^1$
does not intersect the interior of $\cone$.

If two of the weights $\gamma_j^F$ are parallel, then we can
remove one of them without changing the polarized cone $\cone_F$.
Since the first $k'$ weights point in exactly  $k$  different
directions (there are precisely $k$ faces $\sigma'$ of $\cone$
with $\sigma\subset\bar \sigma'$ and $\dim \sigma'=\dim\sigma
+1$), we can assume without loss of generality that $k'=k$.

Let the first $k$ weights span a cone $\tilde C_F$ centered at the origin.
It is clear that  in a small neighborhood of $p$ the polytope
$\cone\cap W$ is equal to the cone $p+\tilde C_F$. Moreover, since the
first $k$ weights are linearly independent, the cones $\tilde C_F$ and
$C_F^1$ are either the same or do not have common points in the
interior. So it remains to show that $\tilde C_F$ and $C_F^1$ are
different, which is the same as showing that during polarization
at least one of the first $k$ weights changes sign.

So, we must show that it is impossible to
have $\gamma_j^F(\xi)<0$ for $1\leq j\leq k$ for every $\xi\in
\bar\Lambda$. If this happens then
the maximum of $q(\xi)$ for $q\in W\cap \cone$ is attained at
$\mu_\cone(F)$, which is impossible, by property (2) of the
cone $\cone$.
\end{proof}

\subsection{Jeffrey-Kirwan localization} Let cones $\Lambda,\bar\Lambda$
and $\cone$ be defined as in the previous section. Consider the symplectic
cut $M_\cone$.
Let $p\in\cone$ be a point  close to the origin. Then
$\mu_\varepsilon=\mu_\cone - \varepsilon p$ for $\varepsilon>0$ is a
moment map on $M_\cone$. Define $\tilde \omega_\varepsilon=\omega_\cone
+i\mu_\varepsilon$ to be an
equivariant symplectic form on $M_\cone$, where $\omega_\cone$ is the
symplectic form on $M_\cone$.

For a form $\eta\in H^*_T(M)$, there
corresponds a form $\eta_\cone\in H^*_T(M_\cone)$. Apply the ABBV
localization theorem to get
\begin{equation}
\label{eq:ABBV-JK}
\int_{M_\cone}\eta_\cone e^{\tilde \omega_\epsilon}=
I^\varepsilon_{M_0}+\sum I^\varepsilon_{F_i} +\sum
I^\varepsilon_{F'_i}
\end{equation}
where
$$
I^\varepsilon_{F} =\frac{1}{d_F}e^{i(\mu_\varepsilon(F))}
\int_F \frac {\iota_F^*(\eta_\cone e^{\omega})}{e(\nu(F))}=
\frac{1}{d_F}e^{i(\mu_\cone(F) -\varepsilon p)}
\int_F \frac {\iota_F^*(\eta_\cone e^{\omega})}{e(\nu(F))}.
$$

Let us now apply Lemma~\ref{lem:minus} to the function $\int_{M_{\cone}}
e^{\tilde \omega_\varepsilon}$ and cones $\bar \Lambda$ and $-\bar\Lambda$:
$$
\res^{\bar \Lambda}([dX] \int_{M_\cone}\eta_\cone e^{\tilde
\omega_\varepsilon})= \res^{-\bar \Lambda}([dX]
\int_{M_\cone}\eta_\cone e^{\tilde \omega_\varepsilon}).
$$
Hence by~(\ref{eq:ABBV-JK}) we get
\begin{equation}
\label{eq:six}
\res^{\bar \Lambda}([dX](I^\varepsilon_{M_0}+\sum I^\varepsilon_{F_i}
+\sum I^\varepsilon_{F'_i} )) =\res^{-\bar
\Lambda}([dX](I^\varepsilon_{M_0}+\sum I^\varepsilon_{F_i}
+\sum I^\varepsilon_{F'_i} ))
\end{equation}
Let us show that four out six terms of~(\ref{eq:six}) are zeros.

Indeed by Lemma~\ref{lem:polarize} the cones $C_{F_i'}$ do not
contain $\varepsilon p$. Hence by property~(\ref{cond1}) of the residue,
we get
$$
\res^{\bar \Lambda}(I^\varepsilon_{F'_i})=0.
$$
Analogously, since the cone $C_{M_0}$ does not contain
$\varepsilon p$
$$
\res^{\bar \Lambda}(I^\varepsilon_{M_0})=0.
$$
Similarly,  for small enough $\varepsilon$
the cones  $-C_{F'_i}$
and $-C_{F'_j}$ do not contain $\varepsilon p$. Hence
by property~(\ref{cond1}) of the residue we have
$$
\res^{-\bar \Lambda}(I^\varepsilon_{F'_i})=0,
\ \ \res^{-\bar \Lambda}(I^\varepsilon_{F_j})=0
$$
for all $i$ and $j$.

Hence only two terms of~(\ref{eq:six}) are not zero. Moreover, since $\bar
\Lambda$ and $\Lambda$ define the same polarization at each $F_i$ one
of these terms can be modified using
$$
\res^{\bar
\Lambda}(I^\varepsilon_{F_i})=\res^\Lambda(I^\varepsilon_{F_i}).
$$

Hence these computations transform~(\ref{eq:six}) into
\begin{equation}
\label{eq:prelimit}
\res^{-\bar \Lambda}(I^\varepsilon_{M_0})=\sum \res^{
\Lambda}(I^\varepsilon_{F_i}).
\end{equation}
Let us now take the limit of both sides of~(\ref{eq:prelimit}) as
$\varepsilon\to 0$. By property (3) of the residue map, and by
equations~(\ref{eq:prop4}) and~(\ref{eq:Euler2}) we have:
$$
\lim_{\varepsilon \to 0^+} \res^{-\bar
\Lambda}(I^\varepsilon_{M_0})= c' \int_{M_0}\kappa_0(\eta
e^\omega),
$$
for some constant $c'$. Hence the limit of~(\ref{eq:prelimit})
gives
$$
\int_{M_0}\kappa_0(\eta e^\omega)= c \sum_i \res^{
\Lambda}\Big(e^{i\mu_\cone(F_i)} \int_{F_i} \frac
{\iota_{F_i}^*(\eta_\cone e^{\omega})}{e(\nu({F_i}))}\Big)
$$
for some constant $c$.

To finish the proof of Theorem~\ref{thm:main}, remember that $F_i$
are the old connected components of the  fixed point set
$M^T_\cone$, so that $F_i\subset M^T$ and
$\mu(F_i)=\mu_{\cone}(F_i)$. In particular,
$\iota^*_{F_i}(\eta_\cone e^{\omega})$ is the same as the
restriction of $\eta e^{\omega}\in H^*_T(M)$ to $F_i$. Moreover,
if for a connected component $F$ of $M^T$ its moment map image
$\mu(F)$ is not inside $\cone$, then the 
cone at $F$ polarized with respect to
$\Lambda$ does not contain the origin, and by
property~(\ref{cond1}) of the residue the term which corresponds
to $F$ in the formula~(\ref{eq:main}) is zero.

%%%%%%%%%%%%%%%%%%%%%%%%%%%%%%%%%%%%%%%%%%%%%%%%%%%%%%%%%%%%%%%%%%%%%%%%
\end{document}